\documentstyle[11pt]{article}
\newtheorem{theorem}{Theorem}[section]
\newtheorem{definition}[theorem]{Definition}
\newtheorem{lemma}[theorem]{Lemma}
\newtheorem{proposition}[theorem]{Proposition}
\newtheorem{corollary}[theorem]{Corollary}
\newenvironment{proof}{\normalsize {\sc Proof}:}{{\hfill $\Box$}}
\def \Z {{\bf Z}}
\def \ZP {{\mathbf Z}^+}
\def \N {{\bf N}}
\def \R {{\bf R}}
\def \F {{\cal F}}
\def \O {{\cal O}}
\def \rtimes {\mbox{$>\!\!\!\lhd$}}
\def \ip {{i'}}
\def \jp {{j'}}
\def \kp {{k'}}
\def \lp {{l'}}

\def \text {\,\,\hbox}

\newcommand{\inv}{^{-1}}

\newcommand{\abs}[1]{\vert #1\vert} 
\newcommand{\e}{\epsilon}
\newcommand{\T}{\mathcal T}
\newcommand{\ncm}[2]{{\left(\!\!\!\begin{array}{c}#1\\#2\end{array}\!\!\!\right)}}

\title{Combing nilpotent and polycyclic groups}
\author{Robert H. Gilman,
Department of Mathematics,\\
Stevens Institute of Technology,
Hoboken, NJ 07030, USA,\\
e-mail: rgilman@stevens-tech.edu\\
\,\\
Derek F. Holt,
Mathematics Institute,\\
University of Warwick,
Coventry CV4 7AL, UK,\\
e-mail: dfh@maths.warwick.ac.uk\\
\,\\
Sarah Rees,
Department of Mathematics,\\
University of Newcastle,
Newcastle NE1 7RU, UK,\\
e-mail: Sarah.Rees@ncl.ac.uk}
\date{}
\begin{document}
\maketitle
\begin{abstract}
The notable exclusions from the family of automatic groups are
those nilpotent groups which are not virtually abelian,
and the fundamental groups of compact $3$-manifolds based on the $Nil$ or $Sol$
geometries.
Of these, the $3$-manifold groups have been shown by Bridson and Gilman to lie
in a family of groups defined by conditions slightly more general than those 
of automatic groups, that is, to have combings which lie in the formal
language class of indexed languages.
In fact, the combings constructed by Bridson and Gilman for these groups
can also be seen to be real-time languages (that is, 
recognised by real-time Turing machines). 

This article investigates the situation for nilpotent  and polycyclic groups.
It is shown that a finitely generated class 2 nilpotent group with cyclic
commutator subgroup is real-time combable, as are also all
2 or 3-generated class 2 nilpotent groups, and groups in specific
families of nilpotent groups (the finitely generated Heisenberg groups,
groups of unipotent matrices over $\Z$ and the free class 2 nilpotent groups).
Further it is shown that any polycyclic-by-finite group embeds in a 
real-time combable group.
All the combings constructed in the article are boundedly asynchronous, 
and those for nilpotent-by-finite groups have 
polynomially bounded length functions, of degree equal to the nilpotency class,
$c$; this verifies a polynomial upper bound on the Dehn functions of those groups
of degree $c$+1.
\end{abstract}
AMS subject classifications: 20F10, 20-04, 68Q40, secondary: 03D40
\renewcommand{\thefootnote}{}
\footnote{
The first author was partially supported by the National Science Foundation.
The second and third authors would like to thank the Fakult\"at f\"ur 
Mathematik of the Universit\"at Bielefeld
for its hospitality while this work was carried out, and the
Deutscher Akademischer Austauschdienst and EPSRC for financial support.}

\section{Introduction}
The aim of this article is to investigate which finitely generated
nilpotent and polycyclic groups have real-time combings, or rather
asynchronous combings which are real-time languages.  

The concept of a combing for a finitely generated group has grown out
of the definition of an automatic group (as introduced in \cite{ECHLPT}).
A formal definition is given in Section~2; informally, a combing
is an orderly set of strands through the Cayley graph of the group,
alternatively a set of words mapping onto the
group for which words
which represent closely related group elements are also
closely related as words; automatic groups possess (synchronous) combings which
are regular languages.
We remark that, by \cite{Bridson2}, the existence of any combing,
synchronous or asynchronous, for $G$ 
implies that $G$ is finitely 
presented, with exponential isoperimetric inequality, and hence soluble word 
problem.

Our interest in combings in other formal language classes arises out of work in
\cite{Bridson&Gilman}, which shows that the fundamental
group of any compact geometrisable 3-manifold has an indexed combing,
that is, a combing
in the formal language class of indexed  languages.
(The indexed languages lie above
the context-free languages in the formal language hierarchy,
and are defined by automata with an attached system of nested stacks.)
It seems natural to ask whether or not the family of
indexed combable groups is large enough
also to contain all finitely generated nilpotent groups.
This paper arose out of an attempt to answer that question.
However, our attention was diverted to the class of real-time languages
by our realisation that all the combings constructed in \cite{Bridson&Gilman}
lay in that class. Further,  while the combing for $\Z^2$ which was
fundamental to the construction of the others in \cite{Bridson&Gilman}
is an indexed language, 
it seems unlikely that the analogous combing for $\Z^n$ is indexed 
(we have not yet proved this), but it is a real-time language. 
The real-time languages, recognised by real-time Turing machines,
are, in particular, recognisable in linear time.
Hence in this paper we address the following question: `Does every finitely
generated nilpotent group have a real-time combing?'

We certainly need to go some way up the language hierarchy to find
combings for nilpotent groups.
By \cite{ECHLPT}, nilpotent
groups, unless virtually abelian, cannot have regular combings, and
by \cite{Bridson&Gilman}, if such a group has a context-free combing, it
cannot be bijective (that is, it cannot contain unique representatives
of each group element). Further, Burillo has proved in \cite{Burillo} that
neither the groups $U_n(\Z)$ of $n$-dimensional unipotent upper-triangular
matrices over $\Z$, nor the $2n+1$-dimensional Heisenberg groups,
defined by the presentations
\begin{eqnarray*}
H_{2n+1}&=&\langle x_1,\ldots x_n,y_1,\ldots y_n,z\,\mid [x_i,y_i]=z,\forall i,\\
& & [x_i,x_j]=[y_i,y_j]=[x_i,y_j]=1,\forall
i,j,i\neq j \rangle
\end{eqnarray*}
can have synchronous, quasigeodesic combings. (Thurston had already
proved in \cite{ECHLPT} that the 3-dimensional Heisenberg group does not
satisfy a quadratic
isoperimetric inequality, and hence cannot be synchronously
combable by quasigeodesics; Burillo's result was proved by consideration
of higher dimensional isoperimetric inequalities.)
In fact both $U_n(\Z)$ and $H_{2n+1}$ have asynchronous combings;
see the end of Section~\ref{sectprops}.

We do not have a complete answer to our question, but have some interesting
partial answers.
We prove that a finitely generated class 2 nilpotent group has  a real-time
combing if it has cyclic commutator group, or can be generated by at most 3
of its elements; we do
so by relating such a group to a semidirect product, and applying a 
result of Bridson (\cite{Bridson}). We give an example of a 4-generated
class 2 nilpotent group which cannot be related to a semidirect product
in this way (but so far we have not proved that this group does not have such a 
combing for some other reason). 
Finally we show that any polycyclic-by-finite group (and hence, of course,
any finitely generated nilpotent group) embeds in a real-time combable group.
We verify that all these combings have fairly good geometrical
properties. In particular, those for class $c$ nilpotent groups have length 
functions which are bounded by a polynomial of degree $c$. By 
\cite{Bridson2}, the existence of such a combing implies a polynomial
upper bound of degree $c+1$ on the Dehn function of the group, verifying a
conjecture commonly attributed to Gersten; this suggests that the 
combings are in some sense optimal.
(The existence of a polynomial upper bound for the Dehn function of
a nilpotent group is proved in \cite{Gersten}, and Gersten's conjecture 
is verified  for
various nilpotent groups, including many of those considered here, in 
\cite{Pittet,Baumslag&Miller&Short}.) 

An entirely different approach towards the construction of combings
is outlined by Gromov in \cite{Gromov}, chapter 5, and described in rather
more detail by Pittet in \cite{Pittet}.  Asynchronous combings for {\em homogeneous}
nilpotent groups are constructed out of {\em homotopic combings} of Lie groups
in which they embed. (Essentially a homotopic combing of a Lie group is a set
of continuous rectifiable paths in the underlying manifold satisfying a 
continuous form of the fellow traveller condition of this paper.)
The construction is geometric, and depends on particular properties of
the Lie group; the language theoretic complexity of these combings is not examined.
Gromov comments that the Heisenberg groups and free class $c$ nilpotent groups 
(the quotients of free groups by the $n$-th terms of their lower central series)
are homogeneous. Pittet deduces, from the combings which can be constructed
for the free nilpotent groups and the groups
\[ G_c = \langle a_1,\ldots a_c,t \mid [a_i,a_j]=1,a_i^t=a_ia_{i+1}, i\neq c,
\,a_c^t=a_c\rangle \]
(defined in \cite{Baumslag&Miller&Short}, and also class $c$ nilpotent),
that these groups consequently have polynomial Dehn functions of degree  $c+1$. 
The results of this article produce real-time combings for the groups $G_c$,
the Heisenberg groups and the class 2 free nilpotent groups, but not obviously  for
the higher class free nilpotent groups.

The work of this paper relates also to work of Baumslag, Shapiro and Short in
\cite{Baumslag&Shapiro&Short}. There the class of parallel poly-pushdown
groups is defined, and proved to contain the fundamental groups of all
compact geometrisable 3-manifolds and
every class 2 nilpotent group. Further, every finitely generated 
torsion-free nilpotent group is proved to embed in a parallel poly-pushdown
group.  Hence the results are analogous to ours. 
However the class of parallel poly-pushdown groups
generalises the concept of an automatic group in a rather different
way; the language associated with a parallel poly-pushdown group satifies
a weaker condition (which can be checked by pushdown automata)
than the fellow traveller condition of this paper, and so
is not necessarily a combing.

\section{Definitions and notations}
Let $G$ be a group, with identity element $1$, and
finite generating set $X$.  Without loss of generality, we may
assume that $X$ is inverse closed, that is, contains the inverse
of each of its elements; we shall make this assumption throughout this
paper.  We call a product of elements in $X$ a word over $X$, and denote by
$X^*$ the set of all such words.
Let $\Gamma=\Gamma_{G,X}$ be the Cayley graph for $G$ over $X$, with vertices
corresponding to the elements of $G$, and,  for each $x \in X$,
a directed edge from the vertex $g$ to the vertex $gx$, labelled by $x$.
Let $d_{G,X}$ measure (graph theoretical) distance between vertices
of $\Gamma_{G,X}$.
For words $w,v \in X^*$, we write $w=v$  if $w$ and $v$ are
identical as words, $w=_G v$ if $w$ and $v$ represent the same element of $G$.
We define $l(w)$ to be the length of $w$ as a string,
and $l_G(w)$ to be the length of the shortest word in $X^*$ 
representing the same element of $G$ as $w$,
that is, the {\em geodesic length} of $w$.
It is straightforward to extend $d_{G,X}$ to a
metric on the 1-skeleton of $\Gamma$. Then each word $w$ can be
associated with a path from $1$ labelled
by $w$, and parameterised by $t \in [0,\infty)$,
such that, for $t<l(w)$, the path from $1$ to $w(t)$
has length $t$, and for $t\geq l(w)$, $w(t)=w(l(w))$.

Suppose that $v,w$ are words in $X^*$, and that $K \in \N$.
We say that $v$ and $w$ {\em asynchronously $K$-fellow-travel} if
there is a differentiable function $h:{\bf R}\rightarrow 
{\bf \R}$, mapping $[0,l(v)+1]$ onto $[0,l(w)+1]$ and strictly increasing 
on $[0,l(v)+1]$, with the property that, for all $t>0$,
$d_{G,X}(v(t),w(h(t)))\leq K$. (The only point of adding 1 to $l(v)$ and 
$l(w)$ in this definition is to deal with the cases where one of the two 
words is trivial.) We call $h$ the {\em relative-speed function}.
If $h$ is the identity function, we say that
$v$ and $w$ {\em synchronously $K$-fellow-travel}.
If for some $M$, and all $t$, $h$ satisfies $1/M \leq h'(t) \leq M$,
we say that $v$ and $w$ {\em boundedly asynchronously fellow-travel} 
with bound $M$. 

We define  a {\em language} for $G$ over $X$ to be a set $L$ of words over $X$
which contains at least one representative for each element of $G$;
$L$ is said to be {\em bijective} if it contains exactly one representative 
of each group element.
We call a language $L$ for $G$ an {\em asynchronous combing} (or just 
{\em combing}) 
if for some $K$
the asynchronous $K$-fellow-traveller condition is satisfied by
all pairs of words $v,w \in L$ for which $w=_G vx$ for some
$x \in X \cup \{ 1 \}$; if relevant pairs of words synchronously
fellow-travel, then $L$ is a {\em synchronous combing},
while if relevant pairs of words boundedly asynchronously fellow-travel,
for some global bound $M$, then $L$ is a {\em boundedly asynchronous combing}.

If $L$ is a combing for $G$ over $X$ then, following \cite{Bridson2},
we define the {\em length function} $f:\N \to \N$ for $L$ by the rule that  
$f(n)$ is the maximum length of a word in $L$ 
of geodesic length at most $n$. The language $L$ is a {\em geodesic 
combing} if $f(n)=n$ for all $n$.

A group $G$ is {\em automatic} if it has a synchronous combing $L$ 
which is a regular language (that is, recognised by a finite state automaton, 
see \cite{Hopcroft&Ullman}), and {\em asynchronously automatic} 
if it has a regular asynchronous
combing. (In fact it is proved in \cite{ECHLPT}, Theorem 7.2.4 that 
any group with a regular, asynchronous combing must have a regular
boundedly asnchronous combing.)
Since (see below) the regular languages form a subfamily
of both the indexed and the
real-time languages, any automatic, or even asynchronously automatic group
is clearly both real-time and indexed combable.

When $G$ is automatic, for each $x \in X \cup \{ 1 \}$, the set $L'_x$ of pairs of
words $v,w \in L$ for which $w =_G vx$ can also be interpreted as a
regular set (over an alphabet of ordered pairs from $X$);
in fact the existence of a regular language $L$ for $G$ for which each such $L'_x$ is 
also regular can be taken as an alternative definition for an automatic group. 
Unfortunately, this alternative view only has limited application when we generalise to other
types of combings.
For synchronous combings, we can still construct finite state
automata to recognise the regular sets $L''_x$ of all pairs of synchronously
$K$-fellow travelling words $w,v$
with $w=_Gvx$ (for some $K$). However the sets $L'_x$  are not regular or even
in the same formal language family as $L$; the regularity
of $L'$ in the case of automatic groups depends on particular properties
of the regular languages (basically, their closure as a language family 
under Boolean operations), which do not hold for other language families.
For an asynchronously combable group, the corresponding sets $L''_x$ are no longer regular,
but automata which read their input
asynchronously from two strings can be built to recognise the fellow traveller property.

A systematic analysis of combings of various types (that is,
associated with a range of fellow-travel properties) and lying in various
formal language classes  is given in \cite{Rees}; in this article we
restrict attention to real-time combable groups,
that is to groups with
asynchronous combings which lie in the formal language class of
real-time languages. Further, the combings we construct
will all be seen to be boundedly asynchronous, and will often have
polynomially bounded length functions.

In this paper we assume familiarity with the basics of
formal language theory, such as the definitions of Turing machines
and finite state automata.
An introduction to the subject, directed towards geometric
group theorists, can be found in \cite{Gilman}; 
\cite{Hopcroft&Ullman} is an excellent standard reference.
Below we give a brief description of real-time languages.
A definition of indexed languages can be found in \cite{Aho}, and of the nested
stack automata which define them in \cite{Aho2}, while the indexed grammars
which define them are described in \cite{Bridson&Gilman}.

Real-time languages (see Rabin's paper \cite{Rabin} for a full definition)
are the languages accepted by deterministic real-time Turing
machines.  These have a fixed number of tapes, one of which
is designated as the input tape, and contains just the input string. The other
tapes can be taken to be infinite in both directions. With each transition,
the machine must read one input symbol and move along the input
tape to the next input symbol.
For each of the other tapes, the machine may write a symbol or a blank, and
then either stay still or move one place to the left or right. The
computation halts when it reaches the end of input. Real-time languages
do not have as many nice closure properties as the others considered in this
paper (for example, they are not necessarily closed under concatenation
with regular languages or under homomorphism). Furthermore, it can be difficult
to determine whether a given language is real-time or not. However, they
do represent a very natural model of linear-time computation, and they
seem to be the most appropriate language for many of the combings that arise
in this paper.
All of $\{ a^{n^2}: n\in \ZP \},\{a^{2^n}:n \in \ZP\}$,
$\{ a^nb^{n^2}: n \in \ZP\}$, $\{ a^{n!} : n \in \ZP\}$ and
$\{ (ab^n)^n: \in \ZP \}$ are real time languages. The first three 
are also indexed 
(\cite{Hopcroft&Ullman}), but the latter two are not (\cite{Gilman2,Hayashi}).
It is proved in~\cite{Rosenberg}
that the language
$$\{ a^{n_1}ba^{n_2}b \ldots a^{n_{r-1}}ba^{n_r}c^sa^{n_{r-s+1}} \, |
r,s, n_1, \ldots, n_r \in \ZP \}$$
is deterministic context-free but not real-time

In the formal language hierarchy, the real-time languages form
a subfamily of the context-sensitive languages which contains
all the regular languages, but (as we see from the example above) not
all context-free languages.
The indexed languages lie between the context-sensitive and 
the context-free languages.
Since every regular language is both real-time and indexed, it is
clear that every automatic, or even asynchronously automatic group is both real-time combable and
indexed combable.

The following, rather surprising, result
shows that real-time combings are only really interesting when they
satisfy additional restrictions.

\begin{proposition}
If $G$ has an asynchronous combing which is recursively enumerable, then $G$
has a real-time asynchronous combing.
\end{proposition}
\begin{proof}
Let $L$ be a recursively enumerable combing for $G$ over $X$ accepted by a Turing machine
$M$. We construct a real-time combing $L'$ for $G$ over a larger alphabet
$X \cup \{ e \}$ by replacing each 
$w \in L$ by a word of the form $we^m$, where $e$ is an alphabet symbol 
representing the identity and $m$ is the sum of the length of $w$ and the number
of moves which $M$ needs to accept $w$.
Note that Proposition~\ref{finvar}, which we shall prove later, implies the existence of a further
real-time combing for $G$ over the original generating set $X$.

$L'$ is accepted by a real-time Turing machine which, given an input word $w'$, first
copies each symbol up to the first occurrence of $e$ onto its work tape, and then operates
as $M$ on the contents of that tape, while continuing to read from
the input tape. It accepts $w'$ provided that it reads only $e$'s while operating as $M$, and that $M$ halts in an accept state just as the end of
the input is reached.
\end{proof} 

However, a real-time combing defined as above is in general not
obviously constructible; its length function $f(n)$ (defined earlier
in this section) might not be recursive.
By contrast, the combings constructed in this paper are
much better behaved.  For instance, they will
all be boundedly asynchronous; this, as we shall now show, implies that
the length function is at worst exponential.
In fact, where the groups involved are nilpotent-by-finite, 
we shall see later that the length function is polynomially bounded.

\begin{proposition}
Let $L$ be a boundedly asynchronous combing for $G$. Then the length
function for $L$ is at worst exponential.
\end{proposition}
\begin{proof}
The boundedness of the combing implies that there is an integer $M$
with the property that, if $h:\R \rightarrow \R$ is the relative-speed
function of two words $v$ and $w$ that asynchronously fellow-travel,
then $(l(w)+1)/(l(v)+1) \leq M$.
Let $c$ be the length of a shortest word $w_0 \in L$ satisfying
$w_0 =_G 1$; the above shows that any other representative
of the identity has length less than $(c+1)M$.
Now let $w \in L$ be a word having geodesic length $n>0$.
We shall show by
induction on $n$ that $l(w) < (c+1)M^n$, which will prove the result.
If $n=1$, then $w$ fellow-travels with $w_0$, so $|l(w)| < (c+1)M$.
If $n>1$, then $w =_G ux$, where $u$ has geodesic length $n-1$ and
$x$ is a generator of $G$. By induction, any word in the combing
for $u$ has length at most $(c+1)M^{n-1}$, but such a word fellow-travels with
$w$, which therefore has length at most $(c+1)M^n$.
\end{proof}
\section{A particularly symmetric combing for $\Z^n$}
\label{Ln}
The combings in this paper are all constructed using a particularly
well behaved combing of the free abelian group $\Z^n$, 
due to Martin Bridson;
it generalises a combing from~\cite{Bridson}. The combing is synchronous 
and geodesic, and behaves particularly well under automorphisms of the group,
due to the fact that it is {\em almost-linear}, that is,
there is a constant $K$ depending only on $n$ such that for any $x\in 
\Z^n$ the combing path in $\Z^n$ from the origin to $x$ and the 
corresponding straight line in $\R^n$ are synchronous $K$-fellow-travellers. 
Here we are considering $\Z^n$ as embedded in the usual way in $\R^n$ and 
thinking of a path in $\Z^n$ as a sequence of elements of $\Z^n$ with 
successive elements distance one apart. 

Let $e_1, \ldots, e_n$ be a basis for $\Z^n$. 
Take $\Sigma_n=\{a_1, a_1\inv, \ldots 
, a_n, a_n\inv\}$, and define $\Sigma_n^*$ to be the free monoid over 
$\Sigma_n$. The correspondence $a_i \to e_i, a_i\inv \to -e_i$ induces a 
 monoid homomorphism from $\Sigma_n^*$ to $\Z^n$. Denote by $\overline 
w$ the image of $ w\in \Sigma_n^*$ under this homomorphism. For example 
$\overline \e = 0$ where $\e$ is the empty word and $0$ is the identity 
element of $\Z^n$. We do not distinguish between $w\in \Sigma_n^*$ and 
the corresponding path in $\Z^n$ from the origin to $\overline w$.

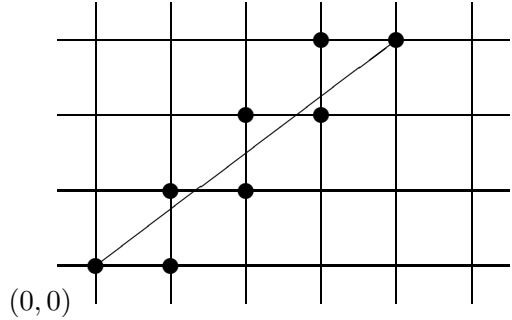
\begin{figure}[htb]
\centering
\bigskip
\unitlength 1mm
\linethickness{0.4pt}
\begin{picture}(65.00,45.00)
\put(10.00,5.00){\line(0,1){40.00}}
\put(20.00,5.00){\line(0,1){40.00}}
\put(30.00,5.00){\line(0,1){40.00}}
\put(40.00,5.00){\line(0,1){40.00}}
\put(50.00,5.00){\line(0,1){40.00}}
\put(60.00,5.00){\line(0,1){40.00}}
\put(5.00,10.00){\line(1,0){60.00}}
\put(5.00,20.00){\line(1,0){60.00}}
\put(5.00,30.00){\line(1,0){60.00}}
\put(5.00,40.00){\line(1,0){60.00}}
\put(10.00,10.00){\line(4,3){40.00}}
\put(10.00,10.00){\circle*{2.00}}
\put(50.00,40.00){\circle*{2.00}}
\put(20.00,10.00){\circle*{2.00}}
\put(20.00,20.00){\circle*{2.00}}
\put(30.00,20.00){\circle*{2.00}}
\put(30.00,30.00){\circle*{2.00}}
\put(40.00,30.00){\circle*{2.00}}
\put(40.00,40.00){\circle*{2.00}}
\put(7.00,7.00){\makebox(0,0)[rt]{$(0,0)$}}
\end{picture}
\caption{The combing path for $(4,3)$ is $a_1a_2a_1a_2a_1a_2a_1$. 
\label{combing}}
\end{figure}

For any point
$p=\sum p_ie_i \in \Z^n$ parameterise the straight line from $0$ to $p$ 
by $x_i=p_ite_i$, $0\le t \le 1$. Start at $t=0$, and each time some $x_i$ 
assumes a positive integer value write down $a_i$, each time some $x_i$
assumes a negative integer value write down $a_i\inv$. If two or more 
$x_i$'s assume integer values at the same time, write down the 
corresponding $a_i$'s in order of decreasing $i$. Figure~\ref{combing} 
shows the combing path for $(4,3)$.

Let $L_n \subset \Sigma_n^*$ be the language defined by the recipe 
above. Clearly $L_n$ is a combing of $\Z^n$, and it is straightforward to 
show that $L_n$ has the desired fellow-traveller property. 

It is shown in \cite{Bridson&Gilman} that $L_2$ is an indexed language.
We devote the remainder of Section~\ref{Ln} to the proof of the 
following result:

\begin{theorem}
\label{realtime}
The language $L_n$ is a real-time language for all $n\geq 1$.
\end{theorem}

\begin{proof}
For $n=1$ the proof is trivial. For $n \geq 2$,
the proof divides into two parts, a reduction to the consideration of
a subset $L_2^+$ of $\Sigma_2^*$,
followed by a demonstration that $L_2^+$ is real-time.

\subsection{Reduction to a 2-dimensional problem}

First we shall reduce to the positive orthant.
For any $S \subset \{1, \ldots, n\}$ define
$f: \Sigma_n^*\to \Sigma_n^*$ to be the monoid homomorphism
which interchanges $a_i$ and $a_i\inv$ for $i \in S$ and fixes all other
$a_i$'s.
A moment's consideration convinces us that $f$ maps $L_n$ to itself. Let
$L_n^+$ be the sublanguage of all paths in $L_n$ which lie in the
positive orthant, i.e., the set of words $w\in L_n$ with $\overline w =
\sum p_je_j$, $p_j \ge 0$. Then $L_n$ is the union of the images of $L_n^+$
under the monoid isomorphisms $f$ corresponding to all $S \subset \{1,
\ldots, n\}$. As real-time languages are closed under union, it follows
that $L_n$ is real-time if $L_n^+$ is.

Next for each $1\le i < j \le n$ let $f_{i,j}:\Sigma_n^* \to \Sigma_2^*$ 
be the monoid homomorphism which sends $a_i$ to $a_1$, $a_j$ to $a_2$, 
and all other generators to $\e$. We claim that $L_n^+$ is the intersection 
of the inverse images of $L_2^+$ under all the $f_{i,j}$'s. Since 
real-time languages are closed under inverse homomorphism and 
intersection, it will follow that $L_n^+$ is real-time if $L_2^+$ is.

It follows from the definition of $L_n$ that $f_{i,j}$ maps $L_n^+$ to 
$L_2^+$. Hence $L_n^+ \subset \cap f_{i,j}\inv(L_2^+)$. For the converse 
suppose $w$ is in the intersection, and $\overline w = \sum p_ke_k$. For 
each $i,j$, $f_{i,j}(w)$ is the combing path for $p_ie_1 + p_je_2$. 
Consequently $f_{i,j}(w)=f_{i,j}(v)$ where $v$ is the combing word in 
$L_n$ for $\sum p_ke_k$. But it is straightforward to check that if 
$f_{i,j}(w)=f_{i,j}(v)$ for all $i,j$, then $w=v$. We conclude that 
$\cap f_{i,j}\inv(L_2^+) \subset L_n^+$.

\subsection{The 2-dimensional problem}

It now suffices to show that $L_2^+$, is real-time. 
From now on, for ease of notation, we shall rewite $a_1$ as $a$ and $a_2$ 
as $b$, and so consider $L_2^+$ as a subset of $\{a,b\}^*$.
Since $b^*$ is a regular language, it is enough to show that
$L_2^\# = L_2^+ \setminus b^*$ is real-time.
We work with $L_2^\#$ because 
it has the following recursive definition.
\begin{eqnarray}
w &=& w_k^n \text{ for some $k\ge 0$ and $n\ge 1$} \label{w} \\
w_j &=& w_{j-1}^{i_j}w_{j-2} \text{ for $2 \le j \le k$ and $i_j \ge 1$} 
\label{recursion}
\end{eqnarray}
and
\begin{eqnarray}
&&\text{ If $k$ is odd, then } w_1 = b^{i_1}a \text{ and } w_0 = b 
\label{evenbasis}\\
&&\text{ If $k$ is even, then } w_1 = a^{i_1}b \text{ and } w_0 = a. 
\label{oddbasis}
\end{eqnarray}
That is, $L_2^\#$ is the collection of all words $w\in \Sigma^*$ satisfying 
some instance of (\ref{w}-\ref{oddbasis}). For example the sequence 
$w_0=a$, $w_1=ab$, $w_2=(ab)^3a$ yields $abababa$, the combing path for 
$(4,3)$.

We shall not verify in this paper that $L_2^\#$ is so defined. This fact 
follows almost immediately from parts of the the proof of Theorem (3.10)
of \cite{Bridson&Gilman}, which verifies that the same language
(known in that paper as $L_1$) is an indexed language.

We construct a real-time Turing machine $\T$ which accepts the language
$L_2^\#$, working with the recursive definition of $L_2^\#$ above.
In our proof, we shall consider 
only inputs beginning with $a$, but our arguments will always extend to 
the other case.

The valid inputs of $\T$ in $a^*b^*$ are  $a$ and $a^{i_1}b$, and 
by~(\ref{recursion}) any valid input not in $a^*b^*$ (and beginning with 
$a$) must have the form $a^{i_1}ba\ldots$. It is straightforward to cope 
with inputs in $a^*b^*$ and to arrange things so that after reading 
$a^{i_1}ba$, $\T$ has one work tape containing $w_0=a$, a second 
containing $w_1=a^{i_1}b$, and so that the tape heads for these tapes 
are positioned at the ends of their tape contents. We shall call 
configurations like this distinguished.

\begin{definition} \label{distinguished} $\T$ is in a distinguished 
configuration if its input is a solution to ~(\ref{w}-\ref{oddbasis}), 
and for some $j \ge 1$
\begin{enumerate}
\item $\T$ has read $w_jw_{j-1}$ from its input;
\item $\T$ has one work tape containing $w_{j-1}$, and a second 
containing $w_j$;
\item The tape heads for these tapes are each positioned either at the 
end of their tape contents or at the beginning.
\item On the work tape containing $w_j$, the squares at distance 
$\abs{w_{j-1}}$ from each end of the tape contents are marked.
\end{enumerate}
\end{definition}

Lemma~\ref{parse} below tells us that if $\T$ is in a
distinguished configuration, then 
either there is a longer input prefix $w_{j+1}w_j$, or the total input 
is $w_{j+1}$ or $w_j^n$ (for some $n>1$). 
$\T$ will verify that the input has one of 
these forms, and in the first case it will simultaneously update the 
work tape containing $w_{j-1}$ so that after coming to the end of the input 
prefix $w_{j+1}w_j$ the work tape contains $w_{j+1}$ and $\T$ is in 
another distinguished configuration. In the second case the input will 
be accepted before $\T$ has time to update the work tape. Thus $\T$ will 
accept all words in $L_2^\#$. Of course for input not in $L_2^\#$, $\T$ will 
reach a point where it is not able to do a required verification, and 
it will reject that input.

During its computation $\T$ will need to compare $w_j$ to segments of 
input. It will do this by traversing the work tape containing $w_j$ in 
either direction. For this method to be feasible we need to know that 
$w_j$ is essentially a palindrome. We shall use the 
following notation. For any word $v$ with length $\abs v \ge 2$, 
$\Phi(v)$ is $v$ with its last two letters reversed. Notice that if 
$u=\Phi(v)$, then $v=\Phi(u)$; and $u\Phi(v)=\Phi(uv)$. Also note that 
$\abs {w_j} \ge 2$ if $j \ge 1$.
Lemma~\ref{Phi}, which refers to any solution
to equations~(\ref{w}-\ref{oddbasis}),
is easy to prove by induction, and we omit the proof.

\begin{lemma}\label{Phi}
For all $j \ge 1$, $w_jw_{j-1} = w_{j-1}\Phi(w_j)$ and
$w_{j-1}w_j = \Phi(w_jw_{j-1})$.
\end{lemma}

To complete our proof of Theorem ~\ref{realtime}, we need Lemmas
\ref{palindrome} and \ref{parse} below, which we also state
without proof. (They follow fairly easily from Lemma~\ref{Phi}.)
They also refer to any solution to 
equations~(\ref{w}-\ref{oddbasis}).

\begin{lemma} For all $j \ge 1$, $w_j$ consists of a palindrome followed 
by $ab$ or $ba$. \label{palindrome}
\end{lemma}

\begin{lemma} \label{parse} Suppose $j>0$, and $w_jw_{j-1}$ is a prefix 
of $w$. One of the following holds.
\begin{enumerate}
\item $j<k$ and $w$ has a prefix $w_{j+1}w_j = 
w_jw_{j-1}\Phi(w_j)^{i_{j+1}-1}w_j$;\label{k-2}
\item $j=k-1, n=1$, and 
$w=w_{j+1}=w_jw_{j-1}\Phi(w_j)^{i_{j+1}-1}$;\label{k-1}
\item $j=k, n>1$, and $w=w_j^n$ consists of $w_jw_{j-1}\Phi(w_j)^{n-1}$ 
with the suffix $w_{j-1}$ deleted.\label{k}
\end{enumerate}
\end{lemma}

Now suppose that $\T$ is in a distinguished configuration as in 
Definition~\ref{distinguished}. $\T$ need only check that its input 
continues in one of the three ways indicated in Lemma~\ref{parse}, 
namely a power of $\Phi(w_j)$ followed by $w_j$, or a power of 
$\Phi(w_j)$, or a power of $\Phi(w_j)$ with the suffix $w_{j-1}$ 
deleted. (Notice that the value of $i_{j+1}$ is determined by the 
occurrence of a substring $w_j$ of the input.) In the first case $\T$ 
must update the worktape containing $w_{j-1}$ so that it contains 
$w_{j+1}$ with the prefix and suffix of length $\abs {w_j}$ marked, and then
position each of its tape heads at one end of the contents of its tape;
afterwards $\T$ is again in a distinguished configuration. In each of
the other two cases $\T$ comes to  the end of its input and accepts the 
input. Actually $\T$ must also verify the parity condition 
of~(\ref{evenbasis}-\ref{oddbasis}), but this task is accomplished by 
checking that the input ends in $a$.

Using Lemma~\ref{palindrome} it is straightforward to design $\T$ so that 
by traversing the worktape containing $w_j$ in either direction it can 
check that the next segment of length $\abs {w_j}$ of the input is 
either $\Phi(w_j)$ or $w_j$. If the input ends with $\Phi(w_j)$, then 
$\T$ accepts. If $w_j$ is encountered, then once we show how the 
appropriate worktape is updated by the time $w_j$ is read from the 
input,  $\T$ will be in another distinguished configuration, namely the 
one corresponding to Definition~\ref{distinguished} with $j$ replaced by 
$j+1$. Finally using the markings on the work tape containing $w_j$, 
$\T$ can also tell when its input ends after a power of $\Phi(w_j)$ with 
the suffix $w_{j-1}$ deleted.

It remains to see how $\T$ updates the worktape containing $w_{j-1}$ 
while it is traversing the worktape containing $w_j$ $i_{j+1}$ times.
Call the first worktape tape 1 and the second tape 2. Suppose that $\T$ 
is at the end of $w_{j-1}$ on tape 1. Since $w_{j+1}= 
w_j^{i_{j+1}}w_{j-1}= w_{j-1}\Phi(w_j)^{i_{j+1}}$, $\T$ need only move 
to the right writing $\Phi(w_j)$ on the tape 1 each time it traverses 
the contents of tape 2. The square on tape 2 which marks the beginning of 
the suffix $w_{j-1}$ of $w_j$ will mark the end of the prefix of length 
$\abs{w_j}$ of $w_{j+1}$ the first time $\Phi(w_j)$ is copied onto tape 
1. In order to mark the appropriate suffix of $w_{j+1}$, $\T$ makes a 
temporary mark at the beginning of each $\Phi(w_j)$ which it writes. The 
first time it scans $w_{j+1}$ in processing the next distinguished 
configuration all except the first of these marks are deleted, and the 
first mark is made permanent. Of course, the marks from the previous iteration
also need to be deleted.

If $\T$ starts at the left of $w_{j-1}$ on tape 1, then the procedure is 
similar. As $w_{j+1}=w_j^{i_{j+1}}w_{j-1}$, $\T$ need only move to the 
left writing $w_j$ from right to left each time it traverses tape 2. As 
$w_j$ is essentially a palindrome, this can be done. The markings are 
done the same way as before with obvious modifications.

This completes the proof of Theorem~\ref{realtime}. 
\end{proof}
\section{Properties of the families of real-time combable
 groups}
\label{sectprops}
We shall construct real-time combings for nilpotent groups by showing that the
groups can be constructed out of free abelian pieces with combings of type 
$L_n$, as already described. 
In this section we prove the necessary closure properties for the families of
real-time combable groups which ensure that this
strategy is valid. In fact, in the following lemma we prove more than we need
(asynchronous results would be enough).
\begin{proposition} 
\label{finvar}
Let $G$ be a finitely generated group.
\begin{description}
\item[(a)] If $G$ is synchronously, asynchronously or boundedly asynchronously
real-time combable, then it is so with respect to any generating set.
\item[(b)] Let $H$ be a subgroup of finite
index in $G$. Then $G$ is synchronously, asynchronously or boundedly
asynchronously 
real-time combable if and only if $H$ is.
\end{description}
In each of the above cases, wherever the original combing has a polynomially
bounded length function, the new combing has a length 
function which is bounded by a polynomial of the same degree.
\end{proposition}
Similar results for asynchronous and synchronous combings in other languages
families are proved in \cite{Bridson&Gilman} and \cite{Rees}.
Although real-time languages do not possess the necessary properties for
all of those results to apply directly, many of the ideas can be used.

\begin{proof}
To prove (a),
let $L$ be a real-time combing for $G$, over a finite (inverse closed)
generating set $X=\{x_1,\ldots x_m\}$,
and suppose that $Y=\{y_1,\ldots y_n\}$ is a second (inverse closed)
generating set for $G$. 
The natural way to define a combing $L'$ over $Y$ is to find a set of
words $w_1,\ldots w_m$ over $Y$, with $w_i$
equal in $G$ to $x_i$, and then define $L'$ to be the set of words
over $Y$ formed by substituting $w_i$ for each occurrence of $x_i$
in each word of $L$. The language
$L'$ naturally inherits asynchronous fellow-traveller
properties from $L$.

We can ensure that $L'$ is recognisable by a real-time Turing machine 
by requiring that each of the substituting words $w_i$ is distinct
and terminated by a string $y_1y_1^{-1}$, which appears nowhere else,
If all the $w_i$ have the same length $k$, then $L'$ satisfies
synchronous fellow traveller properties; although in general we 
cannot organise this, we can at least arrange that all have length 
either $k$ or $k-1$, and pad out the shorter $w_i$'s by strings 
$y_1^{-1}y_1$ on alternate substitutions.

To prove (b), suppose that $H$ is a subgroup of finite index in $G$,
and let $T$ be a finite transversal for $H$ in $G$, containing the
identity element.

If $H$ has a real-time synchronous or (boundedly) asynchronous combing,
then the language for $G$ formed by concatenating that
combing with the elements of $T$ is certainly a real-time
combing of the same type for $G$. 

Now suppose that $L$ is a real-time synchronous or (possibly boundedly) 
asynchronous
combing for $G$ over an (inverse closed) generating set $X$. Let $T$ be a
finite transversal for $H$ in $G$. An appropriate language for $H$ over 
the Schreier generators for $H$ with respect to $X$ and $T$ can be constructed
using the Reidemeister-Schreier rewriting process; the construction
is described in \cite{Rees}. That the language is real-time is
not hard to verify. (The arguments of \cite{Rees} do not in fact always apply
to real-time languages, but can be modified. For instance, where the arguments
require that the language family under consideration is closed under
GSM-mappings, the closure of the family of real-time languages under 
inverse GSM-mappings can be seen to be sufficient for the proof.)

Finally we need to verify that the various
combings constructed  have polynomially bounded length functions, 
given that the same is true of the original combings.

We observe first that the substitutions corresponding to change of 
generators have the effect of changing the length of a word by at
most a constant factor (the maximum length of the old generators
written as words in the new ones).
Similarly, the geodesic length of the substituted word is smaller than
the geodesic length of the original word by at most another
constant factor (the maximum length of the new generators 
written as words in the old).
Hence changing the generating set changes the length function by
at most a constant factor, which can be expressed in terms of the above
constants and the polynomial degree of the original length function.

For each of the remaining cases, the proofs that the length function
continues to be polynomial degree are very similar. Hence we shall give
details of the proof only in one case. Suppose that
$H$ has finite index in $G$, and associated transversal $T$.
Given a combing $L_H$, with polynomial length function $f$, we want to
show that the combing $L_HT$ for $G$ has polynomial length function.

We suppose that $L_H$ is defined over a
generating set $X$ for $H$, and that $Y$ is the set of Schreier
generators for $H$ associated with $X\cup T$ and $T$.
A geodesic representative $w'$ over $X \cup T$ of a word $w$ in $L_H$
can be rewritten as a word of the same length over $Y$ (using the 
Reidemeister-Schreier rewriting process). Its geodesic length
over $X$ is at most a constant factor longer than this (since both $X$ and $Y$
are generating sets for $H$). So, for words representing elements of $H$,
geodesic lengths over $X \cup T$  are at most a constant
factor less than those over $X$, and hence, for some $c$, $cf(n)$
bounds the length of any word in $L_G$ which represents an element of
$H$ and has geodesic length $n$ over $X \cup T$. It remains, for each $t$,
to consider elements of $L_G$ which represent elements of $Ht$.
For such a word $vt$, of geodesic length $n$ over $X \cup T$, 
consideration of $v$, which must have geodesic length at most $n+1$, 
shows that $vt$ has length at most $cf(n+1) +1$.
The result follows.
\end{proof}
\begin{proposition}
\label{directprod}
A direct product of real-time combable groups is real-time combable.
If the factors have boundedly asynchronous combings, then so does
the direct product.
If the combing of each of the factors has polynomially bounded
length function, then so does the combing constructed for the direct product.
\end{proposition}
\begin{proof}
The concatenation of
combings for the direct factors is clearly an asynchronous combing. 
Since disjoint generating sets 
can be chosen for the two factors, the concatenation of the two
languages is easily seen to be in the same family of languages as the 
original languages
(The concatenation of two real-time languages
over non-disjoint alphabets need not be a real-time language).

Bounded asynchronicity is straightforward to check.
The final statement follows from the fact that geodesic words in the
direct product can always be found which are concatenations of geodesic 
words in each of its factors.
\end{proof}
\begin{proposition}
\label{splitext}
Let $G = \Z^n \rtimes H$ be a split extension of an $n$-generated free 
abelian group and a combable group $H$.
Let $L_H$ be the given combing for $H$, and let $L_n$ be the combing for $\Z^n$
described in Section \ref{Ln}.

Then $L_G=L_HL_n$ is a combing for $G$.
If $L_H$ is boundedly asynchronous then so is $L_G$.
If $L_H$ is a real-time language, then so is $L_G$.
\end{proposition}
\begin{proof}
The proof that $L_G$ is a combing is given in \cite{Bridson} (Theorem B).
Since the concatenation of real-time languages over disjoint alphabet sets
is easily seen to be a real-time language, the fact that $L_G$ is real-time
follows immediately from Theorem \ref{realtime}.
That the asynchronicity of the combing is bounded is not proved 
in \cite{Bridson}, but is clear from examination of the proof.
\end{proof}

In order to get a polynomial bound on the length functon for a combing
of a nilpotent group of the form $\Z^n\rtimes H$, we need to look
more closely at the action associated with the extension.
We say that a group $H$ acting on a group $N$  {\em acts
nilpotently} if a series $N=N_0 \supseteq N_1 \supseteq \ldots N_k = 1$ 
of subgroups of $N$ can be found with 
$[N_i,H] \subseteq N_{i+1}$ for $0 \leq i < k$. 
We call such a series an {\em $H$-central series} for $N$. We define the
{\em relative class} of the action (or of the associated
split extension) to be the minimum such $k$. If $N$ is $\Z^n$
and $H$ acts nilpotently on $N$,  with relative class $k$,
then it can be shown that an $H$-central series for $N$ of
length $k$ can be found in which the factors $N_i/N_{i+1}$ are all
torsion-free. This will be proved in Section \ref{nilpotent}, Lemma
\ref{relaction}. Note that if $N \rtimes H$ is nilpotent of class $c$,
then $H$ acts nilpotently on $N$ with relative class at most $c$.
 
\begin{proposition}
\label{polylen}
Let $G  = \Z^n \rtimes H$ be as in the previous lemma, and assume in addition
that $H$ acts nilpotently on $\Z^n$, with relative class $k$, and that
the combing $L_H$ of $H$ has polynomially bounded length function of degree $m$.

Then a generating set for $\Z^n$ can be found such that
the associated combing $L_G=L_HL_n$ of $G$ has length function bounded by a
polynomial of degree at most the maximum of $m$ and $k$.
\end{proposition}
\begin{proof}
Lemma \ref{relaction} ensures the existence of
an $H$-central series
$\Z^n=N_0 \supseteq N_1 \supseteq \ldots N_k = 1$ of $\Z^n$
with torsion-free factors.
Select the generating set $X=\{x_1,\ldots,x_n\}$ for $\Z^n$ 
so that $x_{i_j+1},\ldots,x_n$ generate $N_j$ minimally, for each $j<k$,
$i_0=0$ and $i_k=n$.
Let $c$ be the maximum of the lengths in $X$ of any of the elements $x_j^y$,
where $y$ is an element of the generating set $Y$ of $H$ over
which $L_H$ is defined.
For each $y \in Y$,
if $i_j<i\leq i_{j+1}$, then $x_i^y = x_iw$,
for some word $w$ of length at most $c-1$ in $x_{i_{j+1}+1},\ldots x_n$.
We now establish a bound on the length of $x_i^h$ in terms of $l(h)$, $c$
and $k$.
Let $h=y_{t_1}y_{t_2}\ldots y_{t_l}$ with $y_{t_j} \in Y$. 
Define elements $g_0,g_1,\ldots g_l \in \Z^n$ by $g_0 = x_i$, $g_j = g_{j-1}^{y_{t_j}}$ for $1 \leq j \leq l$, and hence $g_l = x_i^h$.
It is convenient to assign {\em dates} to the occurrences of the generators 
$x_r$ in the words $g_0,\ldots g_l$.
The single occurrence of $x_i$ has date 0, and in general, if an
occurrence of $x_r$ in $g_{j-1}$ has date $m$ and $x_r^{y_{t_j}} = x_rw$, then 
the corresponding occurrences of $x_r$ and the generators of $w$ in $g_j$ have
dates $m$ and $m+1$ respectively.
Let $t(j,m)$ be the total number of occurrences of generators of date $m$ 
in the word $g_j$. Then 
the generators of date $m$ in $g_{j+1}$ include those same $t(j,m)$ generators 
together with at most $(c-1)t(j,m-1)$ new ones arising from the
conjugations by $y_{t_{j+1}}$.
Hence $t(j+1,m) \leq t(j,m) + (c-1)t(j,m-1)$.
From this inequality and the conditions $t(j,0) = 1$ for all $j$
and $t(0,m)=0$ for all $m>0$, it follows by induction that
$t(j,m) \leq \ncm{j}{m}(c-1)^m$ for all $j,m$.
However, since  $t(j,m) = 0$ for all $m\geq k$, we have
$l(x_i^h) \leq a \times l(h)^{k-1}$
for some constant $a$ which depends on $c$ and $k$.

Now let $v$ be a geodesic word over $X \cup Y$, containing $r$
generators from  $X$ and $s$ from $Y$. Let $u$ be the concatenation of
the $s$ generators in $Y$ in the order in which they appear in $v$, and let
$u' \in L_H$  with $u' =_H u$.
Then we have $v =_G u'w'$, where $w'$ is
a product of $r$ elements of $X$ each conjugated by some suffix of $u$.
From the preceding paragraph, there exists $v' \in L_G$
with $v' =_G v$ and with $v'$ equal to $u'$ times a word of
length at most $r\times a \times s^{k-1}$. Since, by assumption,
$l(u')$ is bounded by a polynomial of degree at most $m$ in $l(u)=s$ and
$ras^{k-1} \leq a(r+s)^k$,
the result follows.
\end{proof}

The following corollary follows immediately from repeated application of 
this lemma.
\begin{corollary}
\label{tower}
Suppose that the nilpotent group $G$ is isomorphic to a tower of split
extensions of the form
$\Z^{n_1} \rtimes (\Z^{n_2} \rtimes (\Z^{n_3} \rtimes \ldots \Z^{n_r})\ldots)$.
Then $G$ has a real-time combing $L_G$, whose
length function is bounded by a polynomial with degree the maximum, $k$,
of the relative classes of the extensions.
\end{corollary}
Note that $k$ is no larger than the nilpotency class of $G$.

\begin{corollary}
\label{unipotent}
For any $n$, 
the unipotent groups $U_n(\Z)$ and the Heisenberg groups $H_{2n+1}$,
as defined in Section 1, as well as the free nilpotent groups of
class 2, $Fr_n/\gamma_2(Fr_n)$, are real-time combable. 
The combings are boundedly asynchronous, with length functions
bounded by polynomials of degrees $n-1$, 2 and 2, respectively.
\end{corollary}
\begin{proof}
For the Heisenberg and free nilpotent groups we apply Propositions~\ref{splitext}
and \ref{polylen}; for the unipotent groups we apply Corollary \ref{tower},
observing that $U_n(\Z)$ is a split extension of $\Z^{n-1}$ by $U_{n-1}(\Z)$.
\end{proof}

The following result shows that many soluble groups which 
are far from being nilpotent are also boundedly asynchronously
real-time combable.
\begin{corollary} If $G$ is polycylic, metabelian and torsion-free with
centre disjoint from $G'$, then $G$ has an boundedly asynchronous 
real-time combing.
\end{corollary}
\begin{proof}
That $G$ is polycyclic and metabelian implies that $G'$ is finitely generated.
We now apply a result of Robinson (\cite{Robinson}). The particular form
of this rather general result which we need is stated in \cite{Segal}, namely
that if $A$ is a finitely generated, free abelian normal subgroup of
a group  $G$ such that $G/A$
is finitely generated and nilpotent, and such that $C_A(G)=1$, 
then some subgroup of finite index in $G$ is a split extension of $A$.
We set $A$ to be $G'$, and apply Proposition~\ref{splitext} to get our
result.
\end{proof}.

As an example of a combable group of this form, we have the group
\[ \langle x,y,z \mid  yz = zy, y^x = z,  z^x = yz \rangle \]
which is certainly not automatic (it has exponential isoperimetric 
inequality, see \cite{ECHLPT}, Theorem 8.1.3).

These examples are the building blocks of many others.
The combability of class 2 nilpotent groups with cyclic commutator subgroup,
proved in Section \ref{class2},
is basically a consequence of the combability of the Heisenberg groups;
similarly the embedding theorem for polycyclic-by-finite groups proved in
Section \ref{embedding} follows essentially from the combability of the groups $U_n(\Z)$.

\section{Useful properties of nilpotent groups}
\label{nilpotent}

In this section, we list some facts
about nilpotent groups, which we shall use in the following two sections. 
They are of a standard nature, and only outlines of proofs will be included.
Let
$G=\gamma_1(G) \supset \gamma_2(G) \supset \ldots \supset \gamma_{c+1}(G)
=1$ be the lower central series of $G$, where $G$ is nilpotent of class $c$.
The first result is well-known; see, for example, Lemma~2.6, Corollary~1
of \cite{Hall}, and the second is of a similar nature.
\begin{lemma}
\label{derquot}
If $H$ is a subgroup of the nilpotent group $G$ with $HG'=G$, then
$H=G$.
\end{lemma}

\begin{lemma}
\label{finquot}
Let $H$ be a subgroup of the finitely generated nilpotent group $G$ for which
$|G:HG'|$ is finite.  Then $|G:H|$ is finite.
\end{lemma}
\begin{proof}
Using induction on the class $c$, we may assume that
$|G:H\gamma_c(G)|$ is finite, and so we must show that
$|H\gamma_c(G):H| = |\gamma_c(G):H \cap \gamma_c(G)| <\infty$.

An element of $\gamma_c(G)$
is a product of commutators $[g,k]$ with $g \in G$ and
$k \in \gamma_{c-1}(G)$. Let $|G:HG'| = m$ and
$|\gamma_{c-1}(G):(H \cap \gamma_{c-1}(G))\gamma_c(G)|=l$.
Then $g^m \in HG'$  and $k^l \in (H \cap \gamma_{c-1}(G))\gamma_c(G)$,
and so $[g,k]^{ml} = [g^m,k^l] \in H \cap \gamma_c(G)$. Thus
$\gamma_c(G)/(H \cap \gamma_c(G))$ is finitely generated of finite
exponent, and the result follows.
\end{proof}

The next result allows us to reduce consideration of an arbitrary
nilpotent group $G$ to one where both $G$ and $G/G'$ are torsion
free.
\begin{lemma}
\label{redtotf}
Let $G=\langle g_1,\ldots g_n \rangle$ be a finitely generated nilpotent group.
Then $G$ has an $n$-generator subgroup $G_0$ of finite index such that both
$G_0$ and $G_0/G_0'$ are torsion-free.
\end{lemma}
\begin{proof}
By Theorem 7.8 of \cite{Hall}, there is an integer $r$ such that
$G^r$ is torsion-free. Let $H = \langle
g_1^r, \ldots, g_n^r\rangle$. Then $|G:HG'|$ is finite,
so by Lemma~\ref{finquot}, $|G:H|$ is finite. Since $H \subseteq
G^r$, $H$ is torsion-free.

Now choose $h_1, \ldots, h_m \in H$ such that $h_1H', \ldots h_mH'$
freely generate a maximal free abelian subgroup of $H/H'$, and
let $G_0 = \langle h_1, \ldots, h_m \rangle$. Clearly $m \leq n$.
By Lemma~\ref{finquot}, $|H:G_0|$ and hence $|G:G_0|$ is finite.
Finally, since $G_0' \subseteq H'$, $G_0/G_0'$ is
itself free abelian of rank $m$ and hence torsion-free, so the result follows.
\end{proof}
\begin{lemma}
\label{relaction}
If $H$ acts nilpotently on $N=\Z^n$ with relative class $k$, then 
$N$ has an $H$-central series of length $k$ with torsion-free 
factors.
\end{lemma}
\begin{proof}
Let the given  $H$-central series be
$M_0=1 \subseteq M_1 \subseteq M_2 \subseteq \ldots M_k=N$,
and define the series 
$Z_0=1 \subseteq Z_1 \subseteq Z_2 \subseteq \ldots$ by
letting $Z_{i+1}$ be the inverse image in $N$ of
$C_{N/Z_i}(H)$ for all $i \geq 0$.
Induction shows that $M_i \subseteq Z_i$ for all $i$. Hence $Z_k=N$.

We show that $Z_{i+1}/Z_i$ is torsion-free by induction on $i$.
This is clear for $i=1$.
For $i>1$, 
let $g \in Z_{i+1}\setminus Z_i$ and suppose that $g^m \in Z_i$.
Then there exists $h \in H$ with $[h,g] \not \in Z_{i-1}$;
but $[h,g]^mZ_{i-1}=[h,g^m]Z_{i-1}= Z_{i-1}$,
contradicting the induction hypothesis.
\end{proof}
\section{Class 2 nilpotent groups}
In the following we shall denote by $F_{k,c}$ (or sometimes, where there
is no ambiguity, simply by $F$) the {\em free nilpotent group} on $k$ 
generators of class $c$, that is the quotient of the free group $Fr_k$ on
$k$ generators by its subgroup $\gamma_{c+1}(Fr_k)$.

\label{class2}
Suppose that $G$ is nilpotent of class 2.
In this section we prove that if $G$ has cyclic
commutator subgroup, or can be generated by at most 3 of its elements,
then it is real-time combable.
The method is basically to decompose finite index subgroups of $G$ as split
extensions, and then to apply the results of Section~\ref{sectprops}.
However, we have an example of a 4 generator class 2 nilpotent group
which has no finite index subgroup which is a split extension.

\begin{lemma}
\label{commrels}
Suppose that $G$ is finitely generated nilpotent of class 2, and that $G/G'$
is torsion-free.  Then, for some $k$, 
$G$ can be expressed as a quotient $F_{k,2}/K$ with $K \subseteq F_{k,2}'$.
\end{lemma}
\begin{proof} By Lemma~\ref{derquot}, if $G/G'$ has rank $k$, then $G$
can be generated by $k$ elements. The result follows.
\end{proof}
\begin{theorem}
\label{cyclic}
Any class 2 finitely generated nilpotent group $G$ with $G'$ cyclic has a
real-time combing.  The combing is boundedly asynchronous and has a length
function which is at most quadratic.
\end{theorem}
\begin{proof}
By Lemma~\ref{redtotf} and Proposition~\ref{finvar} (b), 
we may assume that both $G$ and $G/G'$ are torsion free.
Let $G'=\langle c \rangle$.
We claim that some finite index subgroup of $G$ decomposes as a direct product
of a central abelian subgroup $Z$ and a subgroup
$H = \langle a_1,\ldots a_n,b_1,\ldots b_n \rangle$ 
such that for any $i$, $[a_i,b_i]$ is a power of $c$, but for $i \neq j$,
$[a_i,a_j] = [a_i,b_j] = 1.$
Then $H$ decomposes as a 
semidirect product of the subgroups $\langle a_1,\ldots a_n,c\rangle$
and $\langle b_1,\ldots b_n\rangle$, and
the results of Section~\ref{sectprops} imply that $G$ is real-time
combable.
We establish the claim by induction on the
size of a minimal generating set $X$ for $G$.
By Lemma~\ref{derquot}, $|X|$ is equal to the rank of $G/G'$.

First we select any two non-commuting elements
$a_1, b_1 \in X$ and suppose that $[a_1,b_1]=c^i$
For each other $x \in X$,
where $[a_1,x]=c^j$ and $[b_1,x]=c^k$,
the element $y=a_1^kb_1^{-j}x^i$ commutes with $a_1$ and $b_1$.
Let $G_1$ be the group generated by all such elements $y$.
Then $|G:\langle a_1,b_1,G_1\rangle|$ is finite, and
$G_1 \cap \langle a_1,b_1\rangle \subseteq \langle c^i \rangle$.
If $G_1$ is abelian, then $G_1 \cap \langle a_1,b_1\rangle = 1$,
and we have the required decomposition  with $Z=G_1$.
Otherwise we apply induction to $G_1$.

Bounded asynchronicity and the quadratic bound on the length
function follow from Propositions~\ref{finvar} and \ref{splitext},
and Corollary \ref{polylen}.
\end{proof}
\begin{theorem} 
Any two or three generator class 2 nilpotent group $G$ has a real-time combing.
The combing is boundedly asynchronous and has a length function which
is at most quadratic.
\end{theorem}
\begin{proof}
The two generator groups are covered by Theorem \ref{cyclic},
so we assume that $G$ is three generated.
Lemma~\ref{redtotf} and the results of Section~\ref{sectprops} allow us to
assume that $G$ and $G/G'$ are torsion-free.

By Lemma~\ref{commrels}, $G$ is a quotient of $F_{3,2}$ by a normal
subgroup $K$ in the commutator subgroup. Since we may assume that
$G'\cong F'_{3,2}/K$ is non-cyclic and torsion-free, 
and since $F'_{3,2}$ has rank 3,
we need only consider the case where $K$ is 1-generated, that is, where 
$G$ is defined by one relator which is a product of commutators.
In this case, the single relator can be put into the form
$[a,b]^k[a^ib^j,c]^{-1}$. We consider the cases $i=0$ and $i \neq 0$
separately.

First let $i \neq 0$, and let $N= \langle a^ib^j,[a,b]\rangle$,
$H = \langle b,c \rangle $. Then $|G:NH|$ is finite, and
$N$ is abelian and normal in $G$.
Working mod $G'$, we see that, since $a \not \in H$,
$N \cap H \subseteq \langle [a,b]\rangle$.
Hence if $N \cap H \neq 1$ we have a relation between $[b,c]$ and $[a,b]$;
but such a relation cannot be a consequence of the relator
$[a,b]^k[a^ib^j,c]^{-1}$. So $NH$ is a split extension, and hence,
by the results of Section~\ref{sectprops}, $G$ has a real-time combing.

When $i=0$, the one relation can be written as
$[b,a^kc^j]=1$.
Then the group $\langle a^k,b,c^j \rangle$, which has finite index
in $G$, can be written as  a semidirect product of
$N = \langle a^kc^j,[a,c] \rangle$ and
$H = \langle b,a^k \rangle $. Hence again $G$ has a real-time combing.

Bounded asynchronicity and the quadratic bound on the length
function follow again from Propositions~\ref{finvar}, \ref{splitext},
and \ref{polylen}.
\end{proof}

Once we move to four generators, the situation is much less clear.
Let $F$ be the free nilpotent group of class two and rank four, with
generating set $\{a,b,c,d\}$.  Then $F/F'$ is free abelian of
rank 4, and $F' = Z(F)$ is free abelian of rank 6, and is generated by
the six commutators $[a,b]$, $[a,c]$, $[a,d]$, $[b,c]$, $[b,d]$, $[c,d]$.
Let $K$ be the subgroup $\langle [a,b][c,d] \rangle$ of $F$, and
let $G = F/K$. Then $G' = F'/K$ is free abelian of rank 5,
and $G/G' \cong F/F'$.

\begin{proposition}
\label{nodecomp}
No subgroup $G_0$ of finite index in $G$ can be decomposed as a semidirect
product $N \rtimes H$, where $N$ is a nontrivial abelian normal subgroup of
$G_0$.
\end{proposition}

For the proof, assume that there is such a subgroup with a decomposition
of this form. Note that, since $G$ is
torsion-free, $N$ must be a nontrivial free abelian group.
The idea is to show that the free abelian
groups $[N,G_0] \subseteq N$ and $H' \subseteq H$ both have rank 3. 
Since both lie within $G'$, which 
has rank 5, they must then intersect non-trivially. Hence $N$ and $H$ 
intersect non-trivially, which is a contradiction.
This argument follows from a series of lemmas.

\begin{lemma}
\label{comminK}
Let $E$ be any subgroup of $F$ with 3 (or fewer) generators.
Then no nontrivial element of $K$ lies in $E'$.
\end{lemma}
\begin{proof}
Let $E=\langle e,f,g \rangle$
and suppose that $1 \neq ([a,b][c,d])^t \in E'$.
Choose an odd prime $p$ that does not divide $t$, and let
$P = F/F^p.$
Then $P$ is a special $p$-group of order $p^{10}$, with $|P'|=p^6$,
and $P'$ is generated by the six commutators of the
pairs of the generators. Let $Q = EF^p/F^p.$

To simplify notation, we shall
use $a,b,c,d,e,f,g$ to denote the images of these elements in $P$.
Note that $1 \neq ([a,b][c,d])^t \in Q'$.
We can regard $e,f,g$ as elements of the vector space $P/P'$,
and assume that they are in reduced echelon form with respect to the basis
$a,b,c,d$, and we may as well assume that none of $e,f,g$ equals zero
in $P/P'$.  This leaves the following
four essentially different possibilities for $e,f,g$:

$\begin{array}{ll}
(i)\ b, c, d; & (ii)\ ab^i,c,d\;(0 \leq i < p);\\
(iii)\ ac^i,bc^j,d\;(0 \leq i,j <p);& 
(iv)\ ad^i,bd^j,cd^k\;(0 \leq i,j,k<p). 
\end{array}
$

In all cases, $Q' = \langle [e,f], [e,g], [f,g]\rangle$.
Cases (i) and (ii) are impossible, since
none of $[e,f], [e,g], [f,g]$ involves $[a,b]$.
In Case (iv),
we have $[e,f]=[a,b][a,d]^j[b,d]^{-i}, [e,g]=[a,c][a,d]^k[c,d]^{-i}$
and $[f,g]=[b,c][b,d]^k[c,d]^{-j}$. Of these, only $[e,g]$ involves $[a,c]$
and only $[f,g]$ involves $[b,c]$, so $([a,b][c,d])^t$ would have to be a
power of $[e,f]$, which it clearly is not. A similar argument rules out
Case (iii).
\end{proof}
\begin{lemma}
\label{trivcomm}
If $[g,h]=1$ with $g,h \in F$, then
$\langle g,h\rangle F'/F'$ is cyclic.
\end{lemma}
\begin{proof}
Since $F'$ is central in $F$, we may assume that $g=a^ib^jc^kd^l$
and $h=a^\ip b^\jp c^\kp d^\lp$ for some $i,j,k,l,\ip,\jp,\kp,\lp \in \Z$
and then 
$$ [g,h] = [a,b]^{i\jp-j\ip}[a,c]^{i\kp-k\ip}[a,d]^{i\lp-l\ip}
[b,c]^{j\kp-k\jp}[b,d]^{j\lp-l\jp}[c,d]^{k\lp-l\kp},$$
and so we have $i\jp=j\ip$, etc. It can be checked that, for any
solution of these six equations, $g$ and $h$ are powers of a common
element, and the result follows. 
\end{proof}
\begin{lemma}
\label{Nmodrank}
$NG'/G'$ is an infinite cyclic group.
\end{lemma}
\begin{proof} $NG'/G'$ is free abelian, because $G/G'$ is.
We cannot have $N \subseteq G'$, because this would imply
$HG' = G_0$, and then $H' = G_0'$ would have finite
index in $G'$, so $N \cap H$ could not be trivial.
So we may assume that $N$ contains two elements $\bar{g}$ and
$\bar{h}$ such that $\langle \bar{g}, \bar{h} \rangle G'/G'$ is
not cyclic. But then, if $g$ and $h$ are inverse images of  $\bar{g}$
and $\bar{h}$ in $F$, $\langle g,h \rangle F'/F'$ is also
not cyclic, so $[g,h] \neq 1$ by Lemma~\ref{trivcomm}. But $N$ is abelian, so
$[\bar{g},\bar{h}]=1$, which means that $[g,h] \in K$,
contradicting Lemma~\ref{comminK} (applied to $\langle g,h \rangle$).
\end{proof}
\begin{lemma}
\label{Hmodrank}
$HG'/G'$ is free abelian of rank 3.
\end{lemma}
\begin{proof}
Again $HG'/G'$ is free abelian, and
Lemma~\ref{Nmodrank} implies that
its rank is at least 3. But if it had rank 4,
$|G:HG'|$ would be finite.
so (by Lemma~\ref{finquot}) $|G:H|$ and $|G':H'|$ would be finite,
and $N \cap H$ could not be trivial.
\end{proof}
\begin{lemma}
\label{Ncommrank}
$[N,G_0]$ is free abelian of rank 3.
\end{lemma}
\begin{proof}
Let $\hat{N}$ and $F_0$ be the complete inverse images of $N$ and $G_0$ in $F$. 
Then $\hat{N}F'/F'$ is infinite cyclic by Lemma~\ref{Nmodrank};
let $gF'/F'$ be a generator.
Considering the homomorphism $\phi:F/F' \rightarrow [F,\hat{N}]$ defined
by $h \mapsto [g,h]$, we see that $F/C_F(\hat{N})F' \cong [F,\hat{N}]$,
and so these groups are free abelian of the same rank.
By Lemma~\ref{trivcomm}, the group $C_F(\hat{N})F'/F'$ has rank exactly one.
Thus $F/C_F(\hat{N})F'$ 
and $[F,\hat{N}]$ have rank $4-1=3$;
since $|F:F_0|$ is finite, the same is true of $[F_0,\hat{N}]$.
Now the cyclicity of $\hat{N}F'/F'$ implies that
all elements of $[\hat{N},F_0]$ are commutators of the form
$[g^i,h]$.  Then Lemma~\ref{comminK} applied to $E=\langle g,h \rangle$ 
implies that $[\hat{N},F_0] \cap K  = 1$, and the result
follows.
\end{proof}
\begin{lemma}
\label{Hcommrank}
$H'$ is free abelian of rank 3.
\end{lemma}
\begin{proof}
Let $\hat{H}$ be the complete inverse image of $H$ in $F$. Then,
by Lemma~\ref{Hmodrank}, $\hat{H}F'/F'$ is free abelian of rank 3;
let $\hat{H}F'/F' = \langle eF',fF',gF'\rangle$.
Then $\hat{H}' = 
\langle [e,f], [e,g], [f,g] \rangle$, so its rank is at most 3.
If it were less than 3,
then we would have $[e,f]^x[e,g]^y[f,g]^z=1$ for integers $x,y,z$ not
all zero. But it can be shown that this product of commutators is equal to
a single commutator $[h,k]$ for some $h,k \in F$ with $\langle h,k \rangle
F'/F'$ non-cyclic, which contradicts Lemma~\ref{trivcomm}.
To finish, we need only show that $\hat{H}' \cap K = \{ 1 \}$, which
follows from Lemma~\ref{comminK}.
\end{proof}
This completes the proof of theorem \ref{nodecomp}
\section{Embedding in real-time combable groups}
\label{embedding}
In this section we prove the following result.
\begin{theorem}
Any  polycyclic-by-finite group embeds as a subgroup of a real-time
combable group.
Furthermore, any nilpotent-by-finite group embeds as a subgroup of a 
group with a real-time combing with polynomially bounded length function.
\end{theorem}
This result seems interesting in the context of the result of Gersten and
Short (\cite{Gersten&Short}) that no polycyclic
group can embed in a biautomatic group, unless abelian by finite.
\begin{proof}
We start by applying results of Mal\u{c}ev and Hall (\cite{Malcev, Hall}),
which imply that any polycyclic group  (and hence also any
polycyclic-by-finite group) $G$ has a subgroup $G_0$ of finite
index which embeds in a group $T_n(\F)$ of upper triangular matrices
over an algebraic number field $\F$.
In fact, by using a result in an Exercise on page 36 of \cite{Segal},
we can even embed $G_0$ into $T_n(\O)$, where $\O$
is the ring of integers of an algebraic number field.
Note that $T_n(\O)$ is a finitely generated group.
By inducing this representation from
$G_0$ to $G$, we see that $G$ itself embeds in $T_n(\O) wr S_m$, where
$m = |G:G_0|$.
To finish the general, polycyclic-by-finite case,
we need only to observe that $T_n(\O) wr S_m$ is real-time combable.
$T_n(\O)$ can be decomposed as a semidirect product of the form
$E_n(\O) \rtimes (T_{n-1}(\O) \times A)$, where $E_n(\O)$ is
the group of $n \times n$  matrices over $\O$ with 1's on the 
diagonal, and the only other non-zero entries being in the right hand column,
and $A$ is the group of of $n\times n$ diagonal matrices with a unit of $\O$ in
the bottom right hand corner, and all other diagonal entries equal to 1.
Both $E_n(\O)$ and $A$ are clearly finitely generated and abelian; an obvious
induction argument (on $n$), using Propositions~\ref{splitext} and
\ref{directprod}, then proves that $T_n(\O)$ is real-time
combable. Then since $T_n(\O) wr S_m$ contains a direct product of
copies of $T_n(\O)$ as a subgroup of finite index, Propositions~\ref{finvar} and
\ref{directprod} now imply that $T_n(\O) wr S_m$ is real-time combable.
When $G$ is in fact nilpotent-by-finite, $G$ has a subgroup $G_1$,
of finite index $m'$, which is torsion-free nilpotent, and hence
embeds in some $U_{n'}(\Z)$, by \cite{Hall}, Theorem 7.5. 
Then $G$ embeds in $U_{n'}(\Z) wr S_{m'}$,
which is real-time combable with polynomially bounded length
function, by Corollary~\ref{unipotent} and Propositions~\ref{finvar} and
\ref{directprod}.
\end{proof}

\end{document}